\documentclass[12pt,A4Paper]{article}
\title{Catalan equation over $\mathbb{Z}[i]$ for even exponent }
\author{Prem Prakash Pandey, R. Balasubramanian}
\usepackage{amsmath}
\usepackage{amssymb}
\usepackage{amsthm}
\newtheorem{lem}{Lemma}
\newtheorem{rem}{Remark}
\newtheorem{thm}{Theorem}
\newtheorem{prop}{Proposition}
\makeindex
\begin{document}

\maketitle
\begin{abstract}


This article lists all the solutions of the Catalan equation $x^m-y^n=1$ for $x,y \in \mathbb{Z}[i]$, when one of the exponents $m, n$ is even. 

\end{abstract}

\section{Introduction}
The famous conjecture of Catalan (now a theorem due to Mih$\breve{a}$ilescu) states that ``the only solutions of the equation
\begin{equation}\label{1}
x^m-y^n=1 
\end{equation}
with $x,y \in \mathbb{Z},xy \neq 0 \mbox{ and } m,n \in \mathbb{N},m,n>1 \mbox{ are } (\pm 3)^2-2^3=1$".
The conjecture was finally proven by Preda Mih$\breve{a}$ilescu \cite{PM1,PM2,PM3} in 2002. The ingenious solution of Mih$\breve{a}$ilescu was aided by efforts of many authors \cite{LE,VL,KC,JWC1,JWC2,KI}, and rests on his deep insight in theory of  cyclotomic fields. The article of Yuri Bilu \cite{YB1} provides an excellent exposition of the proof. Also the monograph by Rene Schoof \cite{RS} gives a detailed proof together with some historical developments. For more history of the problem we refer the book of Ribenboim\cite{PR}.\\
There are account \cite{PR,BGT} of studies of equation (\ref{1}) over number fields, i.e. finding $x,y \in \mathbb{O}_K,$ the ring of integers of a number field $K$ and rational integers $m,n>1$ satisfying (\ref{1}). From now onwards we will refer equation (\ref{1}) as Catalan equation. The authors in \cite{BGT} showed that over any number field $K$ the Catalan equation has only finitely many solution, as was shown by Tidjeman \cite{RT} for Catalan equation over $\mathbb{Q}$. But the bounds obtained on the possible solutions are astronomical and we are far from listing all the solution of Catalan equation over any number field. The number field analog of the problem has not seen much light of the day. If one wishes to follow the techniques of Mih$\breve{a}$ilescu, then disposing the cases when one of the exponent is even is must. In this article we aim to achieve this when $K=\mathbb{Q}(i)$. This article lists all the non-trivial solutions of the Catalan equation over $\mathbb{Z}[i]$ when one of the exponents is even.\\
\begin{thm}\label{T1}
The only non-trivial solutions to the Catalan equation $$x^m-y^n=1$$ in $\mathbb{Z}[i]$, with $m,n>1$ are  $(-2,3i,3,2),(-2,- 3i,3,2),(3, 2,2,3),(-3, 2,2,3)$,\\
$(i-1,i-2,5,2),(-i-1,-i-2,5,2),(-1-2i,-i+1,2,5),(-1+2i,i+1,2,5),(i-1,2-i,5,2),(-i-1,2+i,5,2),(2i+1,1-i,2,5),(-2i+1,1+i,2,5)$.
\end{thm}
In section 2 we show that to find all solutions to equation (\ref{1}), with $m=p$, a prime bigger than $3$, and $n=2$  it is enough to find all solutions to $x_3^p-4x_2^p=4$. This latter equation is handled in section 3 (for $p \geq 5$). In the course of handling this equation we faced problem of distinguishing $1$ from a primitive $p^{th}$ root of unity and this is achieved in Lemma \ref{L1}. We think Lemma \ref{L1} is a useful result in its own right. In section 4 we solve equation (\ref{1}) for $m=3$ and $n=2$. Then these are combined to list out all the solutions of equation (\ref{1}) in section 5.

\section{Some Reductions}
First we consider the Catalan equation when the exponents are prime, i.e.
\begin{equation}\label{2}
x^p-y^q=1
\end{equation} 
where $x,y \in \mathbb{Z}[i]$ and $p,q$ are primes with $pq$ even. 
If $p=q=2$ then both $x-y \mbox{ and }x+y$ are units and one finds that either $x=0 \mbox{ or }y=0$. Thus we can assume that one of $p \mbox{ and }q$ is even and the other is odd. The equation $x^p-y^2=1$ translates to $y^2-x^p=1$ by the change of coordinates $x\longrightarrow -x \mbox{ and } y\longrightarrow iy$. Thus, its enough to study any one of the equations $x^p-y^2=1 \mbox{ and }x^2-y^q=1$.\\
We consider the equation $x^p-y^2=1$. For $p=3$ this represents an elliptic curve and this will be dealt separately in section 4. So now onwards we will assume $p\geq 5$.\\
%
%
Suppose that $x^p-y^2=1$ has a solution. Then we have
$$x^p=y^2+1=(y+i)(y-i).$$
Claim: $y$ is not a real number.\\
If $y$ is real, then $y$ is a rational integer and $x^p=y^2+1$ is a positive integer. If $y^2+1$ is a $p^{th}$ power of a rational integer then we obtain a solution to the Catalan equation $x^p-y^2=1$ in rational integers, but there are only trivial solutions to this equation over rational integers. If $y^2+1$ is not $p^{th}$ power of a rational integer then the polynomial $X^p-(y^2+1)$ is irreducible over $\mathbb{Z}$ and hence can not have a solution in $\mathbb{Z}[i]$.\\
We consider following two cases \\
Case (1): $y+i \mbox{ and }y-i$ are coprime.\\
This will give $y+i \mbox{ and }y-i$ are $p^{th}$ powers up to a unit. Since all the units in $\mathbb{Z}[i]$ are $p^{th}$ powers, $y+i \mbox{ and }y-i$ are $p^{th}$ powers themselves. 
One has $y+i =x_1^p\mbox{ and }y-i=x_2^p$, which leads to;
\begin{equation}
 x_1^p-x_2^p=2i.
\end{equation}
Since $y$ is not real, $|y+i| \neq |y-i|$ which implies $|x_1| \neq |x_2|$. Without loss of generality we can assume that $|x_1|>|x_2|=\sqrt{n}$ for some positive integer $n$. So one has $|x_1| \geq \sqrt{n+1}$. Now using equation (2) we obtain\\
\begin{align*}
2 = |2i| &= |x_1^p-x_2^p| \\
         & \geq |x_1^p|-|x_2^p| \\
         & \geq (n+1)^{p/2}-n^{p/2} \\
         & \geq 5/2
\end{align*}
which is a contradiction.\\
Case (2): $y+i \mbox{ and }y-i$ are not coprime.\\
Claim: gcd$(y+i,y-i)=2i$\\
A common divisor of $y+i \mbox{ and }y-i$ will divide $2i$. At least one of $y+i \mbox{ and }y-i$ is divisible by $(1+i)^2=2i$, as the power of $1+i \mbox{ in }(y+i)(y-i)$ is at least $p>3$. Also $y+i \mbox{ and }y-i$ differ by $2i$ so the other one too is divisible by $2i$. This proves the claim. \\
Hence one has
$$ y+i=(1+i)^{r_1} x_1^p \quad y-i=(1+i)^{r_2} x_2^p \mbox{ and } x=(1+i)^kx_1x_2, $$
where $r_1,r_2$ are positive integers satisfying min$\{r_1,r_2\}=2$ and $r_1+r_2=kp $ for some positive integer $k$.
Let us assume that min$\{r_1,r_2\}=r_2$, so one gets $(1+i)^{r_1} x_1^p-(1+i)^2 x_2^p=2i$. By putting $x_3=-(1+i)^kx_1$, we obtain $\frac{1}{4}x_3^p-x_2^p=1$ for some integers $x_3,x_2$ in $\mathbb{Z}[i]$. Also $x_3=0$ or $x_2=0$ will lead to $y= \pm i$, which corresponds to a trivial solution. A similar equation unfolds when min$\{r_1,r_2\}=r_1$, namely $4y_1^p-y_3^p=4$, for some $y_1,y_3 \in \mathbb{Z}[i]$. 

\section{The equation $x_3^p-4x_2^p=4$}
Let $\zeta_p$ denote a fixed primitive $p^{th}$ root of unity and $4^{1/p}$ denote the real $p^{th}$ root of $4$.
\begin{lem}\label{L1}
If $x$ and $y$ are non zero integers satisfying $x^p-4y^p=4$ for a prime $p \geq 5$, then the following holds
$$x=4^{1/p}y \sum_{i=0}^{\infty} \binom{1/p}{i} \frac{1}{y^{pi}}.$$
\end{lem}
\begin{proof}
Using binomial expansion we get
 \begin{equation}\label{3}
 x=4^{1/p}y \zeta_p^{n}\sum_{i=0}^{\infty} \binom{1/p}{i} \frac{1}{y^{pi}}.
 \end{equation}
 Claim: $n \equiv 0 \pmod p$.\\
 Assume it is not so, then from equation (\ref{3}) we get
 \begin{equation}\label{4}
 \zeta_p^{-n}=4^{1/p}\frac{y}{x} \sum_{i=0}^{\infty} \binom{1/p}{i} \frac{1}{y^{pi}}.
 \end{equation}
Consider the element $x+4^{1/p}y$, which lies in the ring of integers of $\mathbb{Q}(4^{1/p})$. We wish to show that there is a prime divisor of $x+4^{1/p}y$ which does not divide $px$.\\
If $x+4^{1/p}y$ is a unit, then by taking norm we see that $x^p+4y^p=\pm 1$. This together with $x^p-4y^p=4$ leads to a contradiction. Thus there is a prime dividing $x+4^{1/p}y$. Clearly some prime above $2$ divides $x+4^{1/p}y$. Also it is immediate to see that any prime divisor of $x$ which divides $x+4^{1/p}y$ is a divisor of $2$. Thus it is enough to show that there is a prime divisor of $x+4^{1/p}y$ not dividing $2p$.\\ 
 Assume that the only prime dividing $x+4^{1/p}y$ is a divisor of $2$. By taking norm for the extension $\mathbb{Q}(4^{1/p})/ \mathbb{Q}$ we get $x^p+4y^p=\pm 2^r$, for some positive integer $r$. But comparing the power of $2$ on left we get $r=2$. This does not agree with $x^p-4y^p=4$.\\
Note that $1-\zeta_p$ is totally ramified in $K/ \mathbb{Q}(\zeta_p)$. We will let $\mathbf{p}$ denote the prime ideal of $\mathbb{Q}(4^{1/p})$ above $p$ and $\wp$ will be the prime ideal of $K$ above $p$. Then $1-\zeta_p =\wp ^ p$ and $p=\mathbf{p}^p$.\\
If no prime divisor of $p$ divides $x+4^{1/p}y$ then we are done. Assume that $\mathbf{p} ^t |x+4^{1/p}y$ and the only prime divisors of $x+4^{1/p}y$ are prime divisors of $2p$, then $\wp ^{t(p-1)} | x+4^{1/p} y$ and $\wp ^{t(p-1)} | x+4^{1/p} \zeta_p y$. This implies that $\wp ^{t(p-1)} |1-\zeta_p$ and consequently $t \leq 1$. This forces $t=1$. Now by taking norm in the extension $\mathbb{Q}(4^{1/p}) / \mathbb{Q}$ we get $x^p+4y^p=\pm 4p$. If $x^p+4y^p=4p$ then $x^p-4y^p=4$ gives 
$$x^p(p -1)=4y^p(p + 1).$$
This immediately gives $|y^p| \leq p-1$, from which we arrive a contradiction. Similarly the other case leads to a contradiction.\\
Let $K=\mathbb{Q}(\zeta_p,4^{1/p})$ and $\wp$ be a prime dividing $x+4^{1/p}y$ and not dividing $px$. Let $R$ denote the ring obtained from $\mathbb{O}_K$ by localizing at $\wp$.
From equation (\ref{4}) one sees that $1+\zeta_p^{-n}$ is a power series in $\frac{1}{y^p}$ with coefficient in $R$ and constant term $\frac{x+4^{1/p}y}{x}$. Note that $1+\zeta_p^{-n}$ is a unit in $R$, hence $\frac{x+4^{1/p}y}{x}$ must be invertible in $R$. Consequently $x+4^{1/p}y$ must be invertible in $R$. But since $\wp$ divides $x+4^{1/p}y$, it is not possible. With this contradiction, the claim is established.

\end{proof}
\begin{rem}
The proof of Lemma \ref{L1} is immediate from equation (\ref{4}), by noticing that right side is a real number. But the proof we have given here extends immediately to prove the lemma when $x$ and $y$ are Gaussian integers.
\end{rem}
\begin{thm}\label{T2}
 The only non trivial solution of $x_3^p-4x_2^p=4$, with $p \geq 5$ and $x_3,x_2 \in \mathbb{Z}[i]$, are $p=5,x_3=-(1+i),x_2=i$ and $p=5,x_3=-1+i,x_2=-i$.
 
\end{thm}

\begin{prop}\label{P1}
There are no non trivial solution of $x_3^p-4x_2^p=4$, with $p \geq 7$ and $x_3,x_2 \in \mathbb{Z}[i]$.
\end{prop}
\begin{proof}
Let $x_3,x_2 \in \mathbb{Z}[i]$ be such that $x_3^p-4x_2^p=4$ and $x_3x_2 \neq 0.$\\ If $|x_2|=1,$ then $|x_3|=\sqrt{2}$ and $p=5$. 
Since $p \geq 7$, now onwards we assume that $|x_2|>1$. Hence the binomial expansion 
$$1+\binom{1/p}{1}\frac{1}{x_2^p}+ \ldots $$ of $\large(1+\frac{1}{x_2^p}\large)^{1/p}$ exists. We write 
$$\alpha = 1+\binom{1/p}{1}\frac{1}{x_2^p}+ \ldots .$$ Then from Remark 1 it follows that
$$x_3=x_2 4^{1/p} \alpha .$$
From this we immediately get 
\begin{equation}\label{M}
\left|x_3-4^{1/p}x_2\right| \leq 4^{1/p}\frac{1}{p}\frac{1}{|x_2|^{p-1}} \left|1+\frac{(1/p-1)}{2!}\frac{1}{x_2^{2p-1}}+ \ldots \right|.
\end{equation}
Let us write $x_3=a_3+ib_3$, $x_2=a_2+ib_2$. 
Since $a_3-4^{1/p}a_2=\mbox{Re}(x_3-4^{1/p}x_2)$ so, using inequality (\ref{M}) one obtains 
\begin{equation}\label{M1}
|a_3-4^{1/p}a_2| \leq 4^{1/p}\frac{1}{p}\frac{1}{|x_2|^{p-1}}\left|\left(1+\frac{(1/p-1)}{2!}\frac{1}{x_2^{2p-1}}+ \ldots\right)\right|.
\end{equation}
One knows that $\frac{1}{\sqrt{2}}|x_2| \geq \mbox{min}\{|a_2|, |b_2|\}$. We make the following, \\
Claim: min$\{|a_2|, |b_2|\} \neq 0$. \\
If $a_2 = 0$, then from the inequality (\ref{M1}) we find that $a_3=0$, as the right side quantity in inequality (\ref{M1}) is less than $1$. This will give $$(ib_3)^p-4(ib_2)^p=4,$$ which is not possible, as left hand side is not real.  \\
Considering the imaginary part, from inequality (\ref{M}) we obtain 
\begin{equation}
|b_3-4^{1/p}b_2| \leq 4^{1/p}\frac{1}{p}\frac{1}{|x_2|^{p-1}}\left|\left(1+\frac{(1/p-1)}{2!}\frac{1}{x_2^{2p-1}}+ \ldots \right) \right|.
\end{equation}
If $b_2=0$, then we obtain $b_3=0$. Thus, in this case, $x_2$ and  $x_3$ are real. \\
\\
Define $y$ by $y-i=(1+i)^2x_2^p$, then $y+i=-(1+i)^{-2}x_3^p$ and $y^2+1=(-x_2x_3)^p$. Thus $y$ satisfies $x^p-y^2=1$ with $x=-x_2x_3 \in \mathbb{Z}$. Since $y$ is purely imaginary, by putting $x'=-x$ and $y'=y/i$ we obtain a non trivial integral solution $y'^2-x'^p=1$ of the Catalan's equation over $\mathbb{Z}$ with $p \geq 5$, a contradiction. This contradiction establishes the claim.\\
Let us assume that min$\{|a_2|, |b_2|\}=|a_2|$. Now consider the function $f(x)=x^p-4$, then one has\\
$$\frac{1}{|a_2|^p} \leq |f(\frac{a_3}{a_2})-f(4^{1/p})|=|\frac{a_3}{a_2}-4^{1/p}| |p\xi^{p-1}|,$$ for some point $\xi$ between $\frac{a_3}{a_2} \mbox{ and }4^{1/p}.$
Now using the estimate in (\ref{M1}) we get $$\frac{1}{|a_2|^p} \leq \frac{1}{|a_2|}4^{1/p}\frac{1}{p}\frac{1}{|x_2|^{p-1}}|(1+\frac{(1/p-1)}{2!}\frac{1}{x_2^{2p-1}}+ \ldots)||p\xi^{p-1}|.$$\\ Also we have $|a_2| \leq \frac{1}{\sqrt{2}}|x_2|$, and hence one obtains,\\
$$\frac{1}{|a_2|^p} \leq \frac{1}{|a_2|}4^{1/p}\frac{1}{p}\frac{1}{2^{\frac{p-1}{2}}|a_2|^{p-1}}|(1+\frac{(1/p-1)}{2!}\frac{1}{x_2^{2p-1}}+ \ldots)||p\xi^{p-1}|,$$\\
i.e. $$2^{\frac{p-1}{2}} \leq 4^{1/p} |(1+\frac{(1/p-1)}{2!}\frac{1}{x_2^{2p-1}}+ \ldots)||\xi^{p-1}|.$$\\
Now we note that $1+\frac{(1/p-1)}{2!} \frac{1}{x_2^{2p-1}}+ \ldots$ is dominated by the geometric series $1+\frac{1}{|x_2|^p}+\frac{1}{|x_2|^{2p}}+ \ldots$.\\ 
Using the lower bound on $x_2$, we see that $2^{\frac{p-1}{2}} < 4^{1/p} |\xi^{p-1}|(1+\frac{1}{10}).$ If $|\frac{a_3}{a_2}|<4^{1/p}$ then we obtain $2^{\frac{p-1}{2}} < 4 (1+\frac{1}{10})$ but this is not possible for $p \geq 7$.\\ Now we handle the case $|\frac{a_3}{a_2}|>4^{1/p}$. In this case we have $\xi = 4^{1/p}+\epsilon$, where  $$|\epsilon|  \leq |\frac{a_3-4^{1/p}a_2}{|a_2|}|.$$ So $|\xi|^{p-1} \leq (4^{1/p}+\epsilon)^{p-1} \leq 4^{(p-1)/p}+1.$ To see the last inequality we just notice that 
$|\binom{p-1}{k}4^{(p-1-k)/p}\epsilon^k| \leq 1/p$ and Using this we obtain that $p$ satisfies $2^{\frac{p-1}{2}} \leq 4+4^{1/p}$ and this does not hold for $p \geq 7$.


\end{proof}

\begin{prop}\label{P2}
 Equation $x_3^5-4x_2^5=4$ has no non-trivial solution except $x_3=-(1+i), x_2=i \mbox{ and }x_3= -1+i, x_2=-i$.
\end{prop}

\begin{proof}
We note that if one of $x_3 \mbox{ and }x_2$ is a unit then $x_3=-(1+i), x_2=i \mbox{ or }x_3= -1+i, x_2=-i$. Further, we see that $1+i \nmid x_2$ and hence $|x_2| \geq 2$. Consequently we get $|x_3| \geq 2$.
Note that inequality (\ref{M1}) is valid for $p=5$ too, so we get $|(x_3-4^{1/5}x_2)|<\frac{1}{20}$.\\
We begin with $x_3^5-4x_2^5=4$,\\ i.e. $$(x_3-4^{1/5}x_2)(x_3^4+4^{1/5}x_3^3x_2+\ldots+4^{4/5}x_2^4)=4.$$ \\
Let us put $\frac{4^{1/5}x_2}{x_3}=\eta $, then $|1-\eta|<\frac{1}{40}$ and hence $\frac{39}{40} \leq |\eta| \leq \frac{41}{40}.$ One has 
$$(x_3-4^{1/5}x_2)(x_3^4+ \eta x_3^4+\ldots+ \eta^4 x_3^4)=4,$$
i.e. $$|(x_3-4^{1/5}x_2)|=\frac{4}{5|x_3|^4}\left| \frac{1-\eta}{1-\eta^5} \right|.$$
Similarly one also obtains $$|b_3^5-4b_2^5|=|b_3-4^{1/5}b_2|5|b_3|^4|(1+\tau+\ldots +\tau^4)|$$ for $\tau=\frac{4^{1/5}b_2}{b_3}$. Since $|b_3-4^{1/5}b_2| \leq |(x_3-4^{1/5}x_2)|$, one gets $|1-\tau|<\frac{1}{20}$ and hence $\frac{19}{20} \leq |\tau| \leq \frac{21}{20}.$\\
Using $|b_3-4^{1/5}b_2| \leq |(x_3-4^{1/5}x_2)|$, we have
\begin{equation}
 |b_3^5-4b_2^5| \leq \frac{4}{|x_3|^4}\left|\frac{1-\eta}{1-\eta^5}\right||b_3|^4|(1+\tau+\ldots +\tau^4)|.
\end{equation}
Since $|\tau| \leq \frac{21}{20},$ we have $|1+\tau+\ldots +\tau^4)| \leq 5.6.$\\
Similarly we have $\left|\frac{1-\eta^5}{1-\eta}\right|>4.7$.\\
Now if $|b_3|<0.8|x_3|$, then from the inequality (6) we obtain
$$|b_3^5-4b_2^5| \leq 4(0.80)^4\frac{5.6}{4.7}<2.$$ We show that this is not possible. Since $b_3^5-4b_2^5$ is a non zero rational integer so we need to show that $|b_3^5-4b_2^5| \neq 1$.\\
Let us assume that $|b_3^5-4b_2^5|=1$. Since $x_3^5-4x_2^5=4$, we obtain
$$(a_3+ib_3)^5-4(a_2+ib_2)^5=4.$$ Comparing the imaginary parts, and taking $|b_3^5-4b_2^5|=1$ in account we have\\
$\pm 1 \equiv 0 \pmod 5$, which is not possible.\\
In case $|b_3|\geq 0.8|x_3|$, then $|a_3|<0.6|x_3|$. Now, as done with the imaginary part, using the real part of $x_3-4^{1/5}x_2$ we obtain $|a_3^5-4a_2^5|<1$, this is not possible. This completes the proof of the proposition.
\end{proof}
Theorem \ref{T2} follows from Proposition \ref{P1} and Proposition \ref{P2}.
Proceeding along the same line we obtain
\begin{thm}\label{T3}
 The only non trivial solution of $4y_1^p-y_3^p=4$ with $p \geq 5$ and $y_3,y_1 \in \mathbb{Z}[i]$ are $p=5,y_3=(1+i),y_1=-i$ and $p=5,y_3=1-i,y_1=i$.
\end{thm}

\section{Elliptic curve case}
In this section we settle the equations $x^3-y^2=1 \mbox{ and }x^2-y^3=1$. They both represent elliptic curves defined over $\mathbb{Q}$.\\ 


We will consider the equation $x^2-y^3=1$, which after change of co-ordinate takes the form $y^2=x^3+1$. The first one is dealt similarly. We will let $E$ denote the set of $\mathbb{Q}$-rational points on the curve $y^2=x^3+1$ and $E(i)$ will denote the $\mathbb{Q}(i)$-rational point on the same. Both $E(i) \mbox{ and } E$ have a group structure under `elliptic curve addition +'. Given any point $P=(x,y)$ in $E(i)$, the point $\bar{P}=(\bar{x},\bar{y})$ is also in $E(i)$. Here $z\longmapsto \bar{z}$ is complex conjugation. The point $P+\bar{P}$ of $E(i)$ is stable under complex conjugation and hence is in $E$. Thus we have the trace map $T:E(i)\longrightarrow E$ sending $P \longmapsto P+\bar{P}$. To know the points in $E(i)$ it is enough to find $T^{-1}(P)$ for $P \in E$. Using Cremona's table \cite{JC} we see that that $E$ is of rank $0$ and the torsion group is of order $6$. The six torsion points are $R=(2,3), \mbox{ }2R=(0,1), \mbox{ }3R=(-1,0), \mbox{ }4R=(0,-1), \mbox{ }5R=(2,-3), \mbox{ }6R=(\infty, \
infty)$.\\

Now consider $4R=(0,-1) \in E$, we want to find points $Q \in E(i)$ such that $T(Q)=4R$, i.e. those points $Q=(x,y)$ such that $Q, \bar{Q}$ and $(0,1)$ are collinear. A line $L$ passing through $Q, \bar{Q}$ and $(0,1)$ is given by $y=mx+1$, with $m=\frac{y-\bar{y}}{x-\bar{x}}$. To get the points $Q$ and $\bar{Q}$ we solve the equations $y=mx+1$ and $y^2=x^3+1$. This gives, other than $(0,1)$, a quadratic equation, namely, $x^2-m^2x-2m=0$.\\


Since we are looking for integral points $Q=(x,y)$ and both $Q,\bar{Q}$ lie on $L$ so we have $x+\bar{x}=m$ is an even integer. Further we want the point $Q$ in $E(i)$ and not in $E$ so the equation $x^2-m^2x-2m=0$ shall have two non real roots, also since we want the points to be integral so these roots must be in $\mathbb{Z}[i]$. Hence the discriminant $m^4+8m$ shall be negative of square of an integer. One observes that this is impossible. Thus there are no points on $E(i)$ with $T(P)=4R$. Since $T$ is a homomorphism so there are no points on $E(i)$ whose image under $T$ is $R, 2R$ (if $P\longmapsto R$ then $4P \longmapsto 4R$) and hence also there is no point on $E(i)$ whose image is $5R=-R$.\\
Now consider the case $3R=(-1,0)$. We consider the line through this point, as it is its own reflection, with slope $m$, where $m$ is chosen as in earlier case. The line is given by $y=m(x+1)$, we substitute this in the equation defining the curve to obtain the points of the intersection. We have $(m(x+1))^2=x^3+1$. Canceling the factor$x+1$ we obtain $x^2-(m^2+1)x+(1-m^2)=0$. Again we obtain $m^2+1$ is an even integer and so $m$ is an odd integer. Also the discriminant $(m^2+1)^2-4(1-m^2)$ is negative of square of an integer. This is impossible for any integer $m$. Hence there are no points $P$ on the curve mapping to $3R$ under $T$.\\
Now we consider the last case of point at infinity, the identity of the group law. Here we are looking for points $P$ on the curve $E(i)$ such that $P=-\bar{P}$. If we write $P=(a+ib,k+il)$ then at once we have $b=0,~k=0$. But then from the equation of the elliptic curve we obtain $(il)^2=a^3+1$, i.e. $(-a,l)$ is a solution to $x^3-y^2=1$ in rational integers, this forces $l=0$ and hence $P \in E$. So there are no solution to the equation $x^2-y^3=1$ in $E(i)$ which are not in $E$.\\
For the equation $x^3-y^2=1$ we see that the point at infinity corresponds to one solution $(-2,\pm 3i)$ in $E(i)$. There are no more solution.

\section{Proof of Theorem \ref{T1}}

Because of the reductions of section 2, from Theorem \ref{T2}, Theorem \ref{T3} we see that the equation $x^p-y^2=1$ has no non-trivial solution in $\mathbb{Z}[i]$ for $p \geq 7$. As remarked in section 2, it follows that the equation $x^2-y^q=1$ has no non-trivial solution in $\mathbb{Z}[i]$  for $q \geq 7$. As a consequence, we see that the equations $x^m-y^{2n}=1$ and $x^{2n}-y^m=1$ have no non-trivial solution in $\mathbb{Z}[i]$, whenever the smallest prime divisor of $m$ is at least $7$.\\
In section 3 we found that the only non trivial solution of $x_3^p-4x_2^p=4$ are $-(1+i)^5-4i^5=4$ and $(-1+i)^5-4(-i)^5=4$. Tracing back along the reductions done in section 1, we find non trivial solutions of equation $x^p-y^2=1$ are $(i-1)^5-(i-2)^2=1$ and $(-i-1)^5-(-i-2)^2=1$ respectively. In Theorem \ref{T3} it was shown that the only solution of $4y_1^p-y_3^p=4$ are $4(-i)^5-(1+i)^5=4$ and $4i^5-(1-i)^5=4$. These gives solutions $(-(1-i))^5-(2-i)^2=1$ and $(-1-i)^5-(2+i)^2=1$ respectively. Since both $i-1$ and $-i-1$ are primes in $\mathbb{Z}[i]$ so they can not be of the form $z^t$ for $z \in \mathbb{Z}[i]$ and $t>1$. Also $2+i$ is not of the form $z^t$ for $z \in \mathbb{Z}[i]$ and $t>1$. Consequently the only non trivial solution to equation $x^m-y^{2n}=1,$ when the smallest prime divisor of $m$ is $5$ are $(i-1)^5-(i-2)^2=1,(-i-1)^5-(-i-2)^2=1,(-(1-i))^5-(2-i)^2=1$ and $(-1-i)^5-(2+i)^2=1$. Correspondingly the only solution of $x^{2n}-y^m=1,$ when the least prime divisor of $m$ is $5$, are $(-i+1)^5-(-1-2i)^2=1, (i+1)^5-(1-2i)^2=1,(1+2i)^2-(1-i)^5=1$ and $(-1+2i)^2-(1+i)^5=1$.\\
In section 4 it was shown that the only non-trivial solutions to the equation $x^3-y^2=1 \mbox{ in }\mathbb{Z}[i]$ are $(-2,\pm 3i)$. Hence, the only non-trivial solutions to the equation $x^m-y^{2n}=1$ in $\mathbb{Z}[i]$, when $3|m$ are $(-2,\pm 3i,3,1)$. Similarly we see that $(\pm 3, 2,1,3)$ are the only non-trivial solutions to the equation $x^{2n}-y^m=1$ in $\mathbb{Z}[i]$, whenever $3|m$. Also it was noted in the introduction that the equation $x^m-y^n=1$ has no non-trivial solution in $\mathbb{Z}[i]$, when both $m$ and $n$ are even. 


\noindent{\bf Acknowledgment.}
The first author author is grateful to Professor Joseph Oesterl\'e for suggesting the trace map used in section 4. We will thank Professor R. Thangadurai and R. Mallesham for going through earlier version of draft. We acknowledge the support of research grant from DAE, Government of India.



\end{document}